\documentclass{elsart}

\usepackage{graphics}
\usepackage[all]{xy}
\usepackage{amssymb}

\def\RR{{\mathbb R}}
\def\CC{{\mathbb C}}
\def\NN{{\mathbb N}}

\def\vect{\overrightarrow}
\def\tr{\mbox{\rm Tr}}
\def\bea{\begin{eqnarray}}
\def\eea{\end{eqnarray}}

\def\be{\begin{equation}}
\def\ee{\end{equation}}
\def\og{{``}}
\def\fg{{''}}
\def\Span{\mbox{\rm Span}}
\def\reel{{\mathcal R}{\rm e}\hspace{0.7 mm}}

\begin{document}
\begin{frontmatter}
\title{A Noncommutative View on Topology and Order}
\author{Fabien Besnard}
\ead{Fabien.Besnard@wanadoo.fr}
\address{EPF, 3 bis rue Lakanal, 92330 Sceaux, FRANCE}
\title{}


\author{}

\address{}


\begin{abstract}
In this paper we put forward the definition of particular subsets on a unital $C^*$-algebra, that we call isocones, and which reduce in the commutative case to the set of continuous non-decreasing functions with real values for a partial order relation defined on the spectrum of the algebra, which satisfies a compatibility condition with the topology (complete separateness). We prove that this space/algebra correspondence is a dual equivalence of categories, which generalizes Gelfand-Naimark duality. Thus we can expect that general isocones could serve to define a notion of \og noncommutative ordered spaces\fg. We also explore some basic algebraic constructions involving isocones, and classify those which are defined in $M_2(\CC)$.
\end{abstract}

\end{frontmatter}

\section{Introduction}

Connes' noncommutative geometry is a promising candidate for the merging of General Relativity and the Standard Model of particle theory. However, it faces two major challenges. The first is quantization, and will not be discussed in this paper (see \cite{besnard} and references therein). The second is compatibility with Lorentz signature. One way to deal with this problem would be to perform back and forth Wick rotations. Although this has been proved to be a successful procedure in Minkowski space, it is known to be unreliable at least in some background independent contexts  \cite{loll}. A more satisfying approach is needed that would be an adaptation of the axiomatic description of spectral triples. Several attempts have already been made. A first option is to define semi-riemannian spectral triples, in which one replaces $C^*$-algebras with $*$-algebras of bounded operators on a Krein space \cite{stro}. 
A perhaps less radical possibility is to use the splitting of a globally hyperbolic manifold into space+time. One can then define a collection of noncommutative $3$-spaces with causal relationships between them \cite{hawkins}, \cite{kopf}. It is also possible to arrive at an axiomatic characterization of the collection of local algebras of causal functions on a globally hyperbolic spacetime \cite{moretti}. This last paper has been very inspirational to the author of these lines. Indeed, the axioms that we will state below are already present in \cite{moretti} with the exception of one, which turns out to be crucial,  however.

Thus, the main motivation for this work is to look for a correct setting for \og Lorentzian noncommutative geometry\fg.  However, we will focus only on one aspect of the problem, namely causal order. It is known that all the conformal structure of a spacetime can be recovered from its causal order relation alone \cite{malam}. For this reason, it has long been advocated (see \cite{sork1}) that causal order should be the central object of quantum gravity and that spacetime should be replaced at microscopic scales by what is called a \og causal set\fg, or \og causet\fg\ for short. In the spirit of noncommutative geometry, a dual formulation of a causet in terms of functions preserving the causal order is thus needed. We will arrive at such a formulation in this article. Note, however, that it will be different from the algebraization of causets in terms of indecomposable ideals in the incidence algebra proposed in \cite{sork2}. More precisely, what we seek is a generalization of the Gelfand-Naimark theorem. Just as the Gelfand-Naimark theorem is the point of departure of noncommutative geometry, we can see this generalization to the poset context as the point of departure of noncommutative poset theory.

To be more specific, the kind of space we seek to render noncommutative is a partially ordered set $M$, or poset, endowed with a topology. The two structures have to satisfy a compatibility condition. Several conditions of this kind have been studied. In order to obtain good duality results, the minimal condition to ask is, tautologically, that there are enough continuous non-decreasing functions to fully determine the order relation. Spaces of this kind are called \og completely separated ordered spaces\fg\, but we will more shortly call them \og toposets\fg. We will state some results about toposets in section 2, most of them directly coming from the classic book \cite{nachbin}. Of course, we will need to formulate a set of axioms defining the dual objects. This we will do in section 3, defining the category of $I^*$-algebras. These are couples $(I,C)$ where $C$ is a $C^*$-algebra and $I$ a closed convex cone in $C$, subject to several additional hypotheses inferred from the properties of the cone of non-decreasing functions defined on a toposet. We will then prove our main result in section 4, that is the extension of the Gelfand-Naimark theorem to the category of unital commutative $I^*$-algebras. Unitarity of the algebra means compacity of the underlying space, and one may worry about Lorentzian manifolds being non-compact. But in this case, as we will recall in section 5, one may work with non-decreasing bounded functions, which gives rise to an ordered compactification known as the Nachbin compactification. Section 6 is devoted to the formulation of the definition of noncommutative causal toposets. Finally, in section 7, we will classify the simplest noncommutative $I^*$-algebras one can think of, namely those with $M_2(\CC)$ as algebras.

\section{Completely separated ordered spaces}

First some notations. In a poset, the symbol $\le$ will always denote a partial order relation. We also denote by this symbol the natural order of $\RR$. If in a poset we have two elements $x$ and $y$ such that $x\le y$ and $x\not=y$ we will write $x<y$. If neither $x<y$ nor $y<x$ we write $x\sim y$. 

Let $M$ and $N$ be two posets. We recall that a map $f : M\rightarrow N$ is called isotone iff :
\be
\forall x,y\in M,\ x\leq y\Rightarrow f(x)\leq f(y)\label{isotone}
\ee
It is called antitone if $f(x)\ge f(y)$ in (\ref{isotone}), it is an order-embedding if the implication is an equivalence in (\ref{isotone}), and it is an (order-)isomorphism if it is a surjective order-embedding.

\begin{defn}
A triple $(M,\le,\tau)$ will be called an ordered topological space if $(M,\le)$ is a poset and $(M,\tau)$ a topological space. We set :
$$
I(M,\le,\tau):=\{f : M \rightarrow \RR|f\mbox{ is isotone and continuous}\}
$$
and 
$$I_b(M,\le,\tau):=\{f\in I(M,\le,\tau)|f\mbox{ is bounded  }\}$$

\end{defn}

In the sequel $\tau$ and $\le$ will be understood and we will write $I(M)$ for short. The sets $I(M)$ and $I_b(M)$ may also be defined in case $M$ is only a preordered set. This will be useful later.

The following definition is taken from \cite{nachbin}.

\begin{defn}
A completely separated ordered space is an ordered topological space such that the continuous real isotonies determine the order of $M$. That is to say :
\be
\forall x,y\in M,\ x\le y\Leftrightarrow  \forall f\in I(M),\ f(x)\le f(y)\label{deftoposet}
\ee
\end{defn}

Note in particular that if $f(x)=f(y)$ for all $f\in I(M)$, then $x\le y$ and $y\le x$, thus $x=y$. This shows that the continuous isotone functions separate the points of $M$. However, condition \ref{deftoposet} is strictly stronger than separation of points by $I(M)$. Consider for instance the set $M=[0;1]$ endowed with its usual topology and the partial order $\preceq$ defined by $x\preceq 1/2$ and $1\sim x$ for all $x\in[0;1[$. Then a function $f$ belongs to $I(M)$ if, and only if, it is a continous function which attains its supremum at $x=1/2$. It is clear that $I(M)$ separates the points, but do not determine the order. Indeed, every $f\in I(M)$ satisfies $f(1)\le f(1/2)$ wheras $1\not\preceq 1/2$.

Since the expression \og completely separated ordered space\fg\ is rather lengthy and should have to be used several times, we will call a space of this kind a {\it toposet} (to recall its topological and poset structures).

It is clear that in definition (\ref{deftoposet}) one can use $I_b(M)$ instead of $I(M)$, since whenever $f(x)\le f(y)$ for some $f\in I(M)$, the same relation will hold for the bounded isotone function $\arctan\circ f$.

Obvious examples of toposets are posets with the discrete topology. To see this just consider \og step functions\fg\ defined by $H_x(y)=1$ if $x\le y$, else $H_x(y)=0$. It is readily verified that they are isotone and determine the order.

\underline{Remark} : This example shows that toposet theory is a generalization of poset theory. In fact, it is also a generalization of the theory of functionally Hausdorff spaces\footnote{Sometimes called $T_{2{1\over 2}}$ spaces.}, i.e. spaces in which points are separated by continuous functions. Indeed, a space of this kind is obviously turned into a toposet when endowed with the trivial order relation ($x\le y\Leftrightarrow x=y$).

Another example is given by $(M,\preceq)$ where $M$ is a globally hyperbolic spacetime and $\preceq$ is the causal order relation. This is so because in such spacetimes there exists a global time function $T$ such that $x\prec y\Rightarrow T(x)<T(y)$ \cite{beem}, and if $x\sim y$ it is possible to find global time functions $T_1$ and $T_2$ such that $T_1(x)<T_1(y)$ and $T_2(x)>T_2(y)$. To do so, take $x_\epsilon$ in the future of $x$ and close enough to $x$ so that $x_\epsilon\sim y$ still holds. Then, using the results of \cite{bernalsanchez}, there exists a global time function $T_1$ such that $T_1(x)<T_1(x_\epsilon)=T_1(y)=0$. One finds $T_2$ by exchanging the roles of $x$ and $y$.

Given two posets $(M_1,\le_1)$ and $(M_1,\le_2)$, the product $M_1\times M_2$ can be endowed with the product order $(x,y)\le (x',y')$ iff $x\le_1 y$ and $x'\le_2 y'$. It is easy to see that the product of two compact toposets is a toposet for the product order. On the other hand, the product of two toposets is generally not a toposet for the lexicographic order.

If $M$ is an interval of the real line with the induced topology, one can show that the only total orders which can give $M$ the structure of a toposet are the natural order $\le$ and its dual $\ge$.

Several other compatibility conditions between topology and order have been explored. One is the following :

\begin{defn}
A pospace is an ordered topological space $M$ such that $\{(x;y)\in M\times M|x\le y\}$ is closed in $M\times M$ with the product topology.
\end{defn}

A pospace is sometimes called a $T_2$-ordered space. It is obvious that a pospace is a Haussdorff topological space. It is also obvious that a toposet is a pospace, since the condition (\ref{deftoposet}) determines a closed subset of $M\times M$. It happens that if $M$ is compact, the converse holds (\cite{nachbin}, p 114, corollary). Thus compact toposets and compact pospaces are the same. Note also that they are simply called \og compact ordered space\fg\ by Nachbin.

The next two propositions will be useful later :

\begin{prop}\label{preord}
Let $M$ be a set and $S\subset\RR^M$. Define $\le_S$ on $M\times M$ by $x\le_S y\Leftrightarrow \forall s\in S,\ s(x)\le s(y)$. Then $\le_S$ is a preorder, which is a partial order relation if, and only if, $S$ separates the points of $M$. In this case, if $M$ is a topological space and $S$ only contains continuous functions, then $(M,\le_S)$ is a toposet.
\end{prop}
The proof is immediate. We will call $\le_S$ the (pre-)order defined by $S$. Note that $S$ will not necessarily be equal to $I(M,\le_S)$. A condition for this to happen will be given in corollary \ref{SW2}.

\begin{prop}\label{mred}
Let $(M,\tau,\preceq )$ be a preordered topological set. That is, $\tau$ is a topology, and $\preceq$ a reflexive and transitive relation. Then there exists a topological quotient of $M$, that we will call $M_{red}$,  such that $M_{red}$ is a toposet and the projection $p : M\rightarrow M_{red}$ is isotone. Moreover, $M_{red}$ satisfies the following universal property : for each isotone and continuous $g : M\rightarrow T$, where $T$ is a toposet, there exists $\tilde g : M_{red}\rightarrow T$, isotone and continuous, such that $\tilde g\circ p=g$. Moreover, $M_{red}$ is uniquely determined up to canonical isomorphism by this universal property.
$$\xymatrix{M\ar@{->}[r]^g \ar@{->}[d]_p&T\cr M_{red}\ar@{-->}[ur]_{\tilde g}&}$$
\end{prop}
\begin{pf}
Let us define the equivalence relation $x\equiv y$ iff $\forall f\in I_b(M)$, $f(x)=f(y)$. 

We set $M_{red}=M/\equiv$, equipped with the quotient topology, and call the $p$ the quotient map. Recall that $h : M_{red}\rightarrow X$, where $X$ is any topological space, is continuous if and only if $h\circ p$ is.

For two classes $p(x)$ and $p(y)$ we define :
\be
p(x)\le p(y)\Leftrightarrow \forall f\in I_b(M),\ f(x)\le f(y)
\ee

It is clear that $\le$ is a well-defined partial order relation such that $p$ is isotone.

Now suppose $\psi(p(x))\le \psi(p(y))$ for all $\psi\in I_b(M_{red})$. It is clear by construction that any $f\in I_b(M)$ goes to the quotient and determines $\tilde f\in I_b(M_{red})$ such that $f=\tilde f\circ p$. By hypothesis we thus have $f(x)\le f(y)$ for each $f\in I_b(M)$, and then $p(x)\le p(y)$. This shows that $M_{red}$ is a toposet.

Now let $g : M\rightarrow T$ be as above. Let us show that $g$ goes to the quotient. For this take $x,y\in M$ such that $x\equiv y$. For each $\phi\in I_b(T)$, $\phi\circ g\in I_b(M)$ and thus $\phi(g(x))=\phi(g(y))$. Now since $T$ is a toposet, this implies that $g(x)=g(y)$. Thus $g$ defines $\tilde g$ as stated in the theorem, which is continuous by property of the quotient topology. Now it is easy to show that $\tilde g$ is isotone by the same technique as above. Explicitly :
\bea
&&p(x)\le p(y)\cr
\Rightarrow & \forall f\in I_b(M),\ &f(x)\le f(y)\cr
\Rightarrow& \forall \phi\in I_b(T),& \phi(g(x))\le\phi(g(y)),\mbox{ since }\phi\circ g\in I_b(M)\cr
\Rightarrow& &g(x)\le g(y),\mbox{ since } T\mbox{ is a toposet}\cr
\Rightarrow &&\tilde g(p(x))\le\tilde g(p(y))\nonumber
\eea

There only remains to prove the unicity part. This is usual abstract nonsense, and we leave it to the reader.
\end{pf}

\underline{Remark} : The universal property of $M_{red}$ is equivalent to the following one : for each $\gamma\in I(M)$, there exists $\tilde \gamma\in I(M_{red})$ such that $\tilde \gamma\circ p=\gamma$. Indeed, given that, take a toposet $T$ and $g : M\rightarrow T$ isotone and continuous. Then for any $\phi\in I(T)$, $\phi\circ g:=\gamma\in I(M)$. Then, since $\gamma$ goes to the quotient, it means that $\phi(g(x))=\phi(g(y))$ whenever $p(x)=p(y)$. But this holds for all $\phi$, and $T$ is a toposet. Thus $g(x)=g(y)$ whenever $p(x)=p(y)$.

The theorem below is an adaption of the Stone-Weierstrass theorem to toposets. We will call it the Stone-Nachbin theorem in future references.

\begin{thm}\label{SW}

Let $M$ be a compact toposet, and $J\subset I(M)$ be such that :
\begin{enumerate}
\item $J$ contains the real constants,
\item $J$ is stable by $+$, $\vee$, $\wedge$, multiplication by $\lambda\ge 0$ (i.e. $J$ is a sublattice cone),
\item $J$ determines the order, i.e. $\forall x,y\in M,\ x\le y\Leftrightarrow  \forall f\in J,\ f(x)\le f(y)$
\end{enumerate}
Then $\bar J=I(M)$.
\end{thm}
\begin{pf}
This theorem and its proof can be found in \cite{nachbin} (p. 114, theorem 3) in a slightly different wording. For 
convenience let us repeat the proof here.

If $J$ satisfies the hypotheses of the theorem then $\bar J$ also does. We can then suppose that $J$ is closed and prove that $J=I(M)$. We take  $F$ in $I(M)$. All we need to do is to apply the lattice version of the Stone-Weierstrass theorem, and for this, we have to show that for all $x,y\in M$ with $x\not=y$, there exists an $f_{x,y}$ in $J$ such that $f_{x,y}(x)=F(x)$ and $f_{x,y}(y)=F(y)$. If $F(x)=F(y)$ then we can take a constant. If not, we can suppose without loss of generality that $F(x)<F(y)$. Since $F$ is isotone, this implies that $x\not\ge y$. But $J$ determines the order, so there exists $j\in J$ with $j(x)<j(y)$. Now define :
$$\lambda:={F(y)-F(x)\over j(y)-j(x)};\qquad \mu:=F(x)-{F(y)-F(x)\over j(y)-j(x)}j(x)$$
We have $\lambda>0$. Thus $f_{x,y}:=\lambda j+\mu 1\in J$ and has the required property.
\end{pf}

The following immediate corollary is in fact the initial formulation of theorem 3 in \cite{nachbin}.

\begin{cor}\label{SW2}
Let $M$ be a compact topological space and $J\subset{\mathcal C}(M)$. Then $(M,\le_J)$ is a toposet such that $I(M)=J$ if and only if $J$ is stable by $+$, $\vee$, $\wedge$, multiplication by $\lambda\ge 0$, contains the constants, separates the points, and is closed.
\end{cor}

Finally, there is a generalization of Tietze's extension theorem in compact toposets.

\begin{thm}\label{extension}
Let $M$ be a compact toposet and $F$ a closed subset of $M$. Then every element of $I(F)$ can be extended to an element of $I(M)$.
\end{thm}

This ``Tietze-Nachbin'' theorem is a corollary of \cite{nachbin}, theorem 4 p. 48, which states that compact ordered spaces (i.e. compact toposets) are normally ordered, and theorem 6 p. 49 which formulates the extension property for normally ordered spaces.

\section{Isocones and $I^*$-algebras}
Let $C$ be a unital $C^*$-algebra. We will use the following notations : the set of hermitean elements of $C$ is $\reel C:=\{a\in C|a=a^*\}$ and the set of positive elements of $C$ is $C_+:=\{a\in\reel C|\sigma(a)\subset\RR_+\}$ where $\sigma(a)$ denotes the spectrum of $a$. We will write $a\ge 0$ if $a\in C_+$. It is well known that $a\ge 0$ iff $\exists b\in C,\ a=b^*b$ and that for any $a\ge 0$, there exists a unique $b\ge 0$ such that $a=b^2$, and we write $b=\sqrt{a}$. We also recall that for any $a\in C$ one defines $|a|=\sqrt{a^*a}$.

At times we will need the set $C_+^\times:=\{a\in \reel C|\sigma(a)\subset\  ]0;+\infty[\}$. So $a\in C_+^\times$ iff $a$ is both positive and invertible. It is easy to see that $C_+^\times$ is an open subset of $\reel C$.

Let $a,b\in\reel C$. We define :
$$a\wedge b={a+b\over 2}-{|a-b|\over 2};\qquad a\vee b={a+b\over 2}+{|a-b|\over 2}$$

We must make some comments here. The two operations just defined are not in general meet and join in the sense of the partial order induced by $C_+$. For instance, it is known \cite{kad} that two self-adjoint bounded operators on a Hilbert space have a least upper bound only when they are comparable. One must also pay attention to the fact that $\vee$ and $\wedge$ do not generally satisfy the lattice axioms. 
They do only when the algebra is commutative \cite{sher}. 

In view of the Stone-Nachbin theorem, the following definition is a natural noncommutative generalization of the set of isotone functions on a toposet.
 
\begin{defn}
A (strong) $I^*$-algebra is a couple $(I,C)$ such that :
\begin{enumerate}
\item $C$ is a unital $C^*$-algebra,
\item $I$ is a closed convex cone of $\reel C$, containing the constants and stable by $\wedge$ and $\vee$, i.e. :
\begin{enumerate}
\item $I\subset \reel C$,\label{deuxa}
\item $\overline{I}=I$,\label{ferme}
\item $\forall b,b'\in I$, $b+b'\in I$,\label{somme} 
\item $\forall \lambda\in\RR_+$, $\forall b\in I$, $\lambda b\in I$,\label{scal}
\item $\forall b,b'\in I$, one has $b\vee b'\in I$ and $b\wedge b'\in I$,\label{meetandjoin}
\item $\forall x\in\RR$, $x1\in I$.\label{constants}
\end{enumerate}
\item $\overline{\Span(I)}=C$\label{dense}
\end{enumerate}
\end{defn}

Examples of $I^*$-algebras are given by couples $(I(M),{\mathcal C}(M))$ where $M$ is a compact toposet. The direct part of corollary \ref{SW2} shows that the second set of properties is satisfied. To see that the third property also is consider $V=\Span_\RR(I(M))$. This is a real algebra of continuous functions on a compact space. Indeed, an element $h$ of $V$ may be written $h=f-g$ where $f,g\in I(M)$. Moreover, using the compacity of $M$ we can write $f=f_+-\lambda 1$ and $g=g_+-\mu 1$ where $f_+$ and $g_+$ are positive functions, and $\lambda,\mu\in\RR_+$. Thus we have $h=(f_++\mu 1)-(g_++\lambda 1)$, so any element of $V$ can be expressed as a difference of two positive elements of $I(M)$. This easily entails that $V$ is stable by products. Since it also contains the constants and separates the points by the toposet property, it is dense in $\reel C(M)=C(M,\RR)$ by the algebra version of the Stone-Weierstrass theorem. But $C(M)=\reel C(M)+i\reel C(M)$ and so $\Span_\CC(I(M))$ is dense in $C(M)$.

For an obvious, but trivial, example of a noncommutative $I^*$-algebra one can take $(\reel A, A)$ where $A$ is any noncommutative $C^*$-algebra.

\underline{Remark}: The \og causal cones\fg\ defined in \cite{moretti} satisfy some of the axioms of isocones, but \ref{meetandjoin} is lacking. The reason is that Moretti's setting contains other structures, satisfying their own axioms, which interact with causal cones. In the end, the causal structure of a globally hyperbolic manifold is recovered, among other things. In particular this implies that causal cones are indeed isocones. If we are concerned with the causal structure only, and do not have other structures at our disposal, then the axioms of isocone are necessary and sufficient, as we will see. 

Of course, any condition which reduces to the stability under meet and join when the algebra is commutative could be required instead of \ref{meetandjoin}. In particular we could ask for the less restrictive condition :
$$(2e') : \forall b,b'\in I\mbox{ such that } [b,b']=0, \mbox{ one has } b\vee b'\in I\mbox{ and } b\wedge b'\in I$$

The corresponding objects will be called weak $I^*$-algebras. In the sequel, when we make a statement about $I^*$-algebra (without adjective), it will be understood that it stands for both strong and weak $I^*$-algebras.

\underline{Remark 1} : Let us call $I_+=I\cap C_+$. Then $I=I_++\RR_-1$. Indeed, for all $x\in I$, $x+\|x\| 1\in I_+$.

\underline{Remark 2} : Using the identity $(a+k 1)\wedge (b+k 1)=a\wedge b+k 1$, with $k\in\RR$, and similarly for $\vee$, we see that using axioms \ref{somme} and \ref{constants} we can weaken axiom \ref{meetandjoin} or 2e'   by assuming that $b,b'\in C_+^\times $.

A subset of a $C^*$-algebra  satisfying the axioms \ref{deuxa}  to \ref{constants} above will be called a {\it pre-isocone} (strong or weak, depending on which version of axiom \ref{meetandjoin} is used). If it moreover satisfies the third axiom, it will be called an {\it isocone}. Thus an $I^*$-algebra is a $C^*$-algebra equipped with an isocone. A sub-$I^*$-algebra of $(I,C)$ is a couple $(J,B)$ such that $B$ is a sub-$C^*$-algebra of $C$, $J\subset I$ is a pre-isocone and $(J,B)$ is an $I^*$-algebra.

Let us now see how the usual algebraic constructions extend to the $I^*$-algebra setting.

\begin{defn}\label{istarmorph}
An $I^*$-morphism from $(I,C)$ to $(I',C')$, where the latter are two strong $I^*$-algebras, is a mapping $\rho$ from $C$ to $C'$  such that :
\begin{enumerate}
\item $\rho$ is a $*$-morphism between $C$ and $C'$,
\item $\rho(I)\subset I'$,
\item For all pre-isocone $J\subset I$, $\rho(J)$ is closed.
\end{enumerate}
When $(I,C)$ and $(I',C')$ are two weak $I^*$-algebras, an $I^*$-morphism $\rho$ from $(I,C)$ to $(I',C')$ is required to satisfy 
\begin{enumerate}\setcounter{enumi}{3}
\item For all pre-isocones $J\subset I$, $\forall x,y\in J$, $[x,y]\in\ker\rho\Rightarrow x\wedge y$ and $x\vee y\in J+\ker\rho$.
\end{enumerate}
in addition to the three above conditions. 
\end{defn}

As a result of this definition, if $\rho$ is an $I^*$-isomorphism, one has $\rho(I)=I'$. We also note that requirements 3 and 4 are trivially satisfied when $\rho$ is injective.

\begin{prop}\label{imagemor} Let $\rho$ be a morphism from $(I,C)$ to $(I',C')$. Then the couple $(\rho(I),\rho(C))$ is an $I^*$-algebra. Moreover, for all pre-isocone $J\subset I$, $\rho(J)$ is a pre-isocone of $C'$.
\end{prop}
\begin{pf}
We know from the definition of $I^*$-morphims that $\rho(I)$ is closed. Moreover, it is immediate that $\rho$ automatically respects $\vee$ and $\wedge$ as a $*$-morphism. All verifications are then trivial in the \og strong\fg\  context. For a morphism of weak $I^*$-algebras and a pre-isocone $J\subset I$, let us show that $\rho(J)$ satisfies 2e'. For this take $\rho(x)$ and $\rho(y)$ such that $x,y\in J$ and $[\rho(x),\rho(y)]=\rho([x,y])=0$. Thus, $[x,y]\in \ker\rho$, and the fourth condition ensures that $\rho(x)\wedge\rho(y)\in\rho(J)$, and similarly for $\vee$.
\end{pf}

\begin{prop}
The class of strong (resp. weak) $I^*$-algebras with their morphisms is a category.
\end{prop}
\begin{pf}
Let $\rho : (I,A)\rightarrow (J,B)$ and $\rho' : (J,B)\rightarrow (K,C)$ be two morphisms between $I^*$-algebras. We need to prove that $\rho'\circ\rho$ is an $I^*$-morphism. The two first conditions are obvious. Let $F$ be a pre-isocone of $A$. Then $\rho(F)$ is a pre-isocone of $B$ and lies in $J$. Thus $\rho'(\rho(F))$ is closed, and the third condition is satisfied.

Finally, in the weak case, we consider $x,y\in F$ such that $[x,y]\in \ker\rho'\rho$. Then $[\rho(x),\rho(y)]\in\ker\rho'$. Now since $\rho(F)$ is a pre-isocone of $J$, we have $\rho(x)\wedge\rho(y)=\rho(x\wedge y)\in\rho(F)+\ker\rho'$. Thus, $x\wedge y\in \rho^{-1}(\rho(F)+\ker\rho')=F+\ker\rho'\rho$, and the result follows. The same holds with $\vee$ instead of $\wedge$.
\end{pf}

We will call ${\mathcal IC}$ the category of strong $I^*$-algebras, and ${\mathcal WIC}$ the category of weak ones. 

\underline{Remark} : Take two $I^*$-algebras of the form $(\reel A,A)$, $(\reel B,B)$. Then any $*$-morphism between $A$ and $B$ gives rise to an $I^*$-morphism if and only if the third condition is fulfilled. Were it always the case, strong $I^*$-algebras would be a generalization of $C^*$-algebras in the same way as compact toposets are a generalization of compact sets. Or to put it in categorical terms, $C^*$-algebras and their morphisms would form a full subcategory of ${\mathcal IC}$. We do not know whether this is true or not.

By the Gelfand-Naimark-Segal theorem, every unital $C^*$-algebra is isomorphic to a sub-$C^*$-algebra of the algebra ${\mathcal B}({\mathcal H})$ of bounded operators on some Hilbert space ${\mathcal H}$. By proposition \ref{imagemor} and the remark above that any injective $*$-morphism satisfies points 3 and 4 of the definition of $I^*$-morphisms, we obtain :

\begin{prop}
Any $I^*$-algebra is a sub-$I^*$-algebra of $(\reel {\mathcal B}({\mathcal H}),{\mathcal B}({\mathcal H}))$ for some Hilbert space ${\mathcal H}$.
\end{prop}

We next define ideals so that they are the kernels of morphisms.

\begin{defn}
Let $(I,C)$ be an $I^*$-algebra. A subset ${\mathcal N}$ of $C$ is called an $I^*$-ideal of $(I,C)$ (or simply an ideal, if the context is clear) if it is a closed two-sided ideal of $C$ such that $J+{\mathcal N}$ is closed for every pre-isocone $J\subset I$. If $(I,C)$ is a weak $I^*$-algebra, we additionaly ask that for all pre-isocone $J\subset I$, and for all $x,y\in J$, $[x,y]\in{\mathcal N}\Rightarrow x\wedge y$ and $x\vee y\in J+{\mathcal N}$.
\end{defn}

With this definition and  proposition \ref{imagemor}, the quotient of an $I^*$-algebra by an ideal is an $I^*$-algebra.

It is a simple matter to verify that for two isocones $I$ of $C$ and $I'$ of $C'$, the direct sum $I\oplus I'$ is an isocone of $C\oplus C'$. If $(I,C)=(I(M),C(M))$ and $(I',C')=(I(M'),C(M'))$ then $I\oplus I'=I(M\coprod M')$ where  the order on the disjoint union of the toposets $M$ and $M'$ is the unique order extending the ones on $M$ and $M'$ and such that $M\sim M'$.

If $B$ is a unital sub-$C^*$-algebra of $C$, and $(I,C)$ is an $I^*$-algebra, then $(I\cap B, B)$ is not necessarily an $I^*$-algebra. The only axiom that is not always satisfied\footnote{For a counterexample, take $C={\mathcal C}([0;1])$ and $I$ the isocone of continuous non-decreasing functions for the natural order, and take $B=\{f\in C|f(0)=f(1)\}$.} is $3$.

When it is the case, we will say that $B$ is transverse to $I$. If $C^*(a)$, the $C^*$-algebra generated by an element $a$ of $C$, is transverse to $I$, we will say that $a$ is transverse to $I$.

Let $S$ be a subset of $C_+$ such that $\overline{\Span(S)}=C$. Since $\reel C$ is an isocone containing $S$, the intersection $\langle S\rangle$ of all isocones containing $S$ will be non-empty. It is clearly an isocone, and we will call it the isocone generated by $S$. Let $(I,C)$ and $(I',C')$ be two $I^*$-algebras. Then for any notion of tensor product, $\langle I_+\otimes I'_+\rangle$ defines an isocone in $C\otimes C'$.

For instance, using the Stone-Nachbin theorem, one can prove that if $C={\mathcal C}^0([a;b])$, and $I$ is the set of continuous non-decreasing functions on $[a;b]$, then $\langle I_+\otimes I_+\rangle$ is the set of continous functions on $[a;b]^2$ that are non-decreasing for the product ordering.

\section{The Gelfand-Naimark theorem for $I^*$-algebras}

The class of toposets (resp. compact toposets or equivalently compact pospaces) together with their continuous isotonies forms a category. This category will be called ${\mathcal T}$ (resp. ${\mathcal KT}$).

If $f : M\rightarrow N$ is a continuous isotony between two toposets then we write $f^*$ for the pullback by $f$ from $I(N)$ to $I(M)$. That is : $f^*(g)=g\circ f$.

\begin{thm} Let $(M,\preceq)$ and $(N,\le)$ be two compact toposets and $f : M\rightarrow N$ a continuous isotony. Then write ${\mathcal I}(M)=(I(M),C(M))$ where $C(M)$ is the set of continuous maps from $M$ to $\CC$, and write ${\mathcal I}(f)=f^*$. Then ${\mathcal I}$ is a contravariant functor from ${\mathcal KT}$ to $ab-{\mathcal IC}$.
\end{thm}
\begin{pf}
The only thing we still need to prove is that $f^*$ is a morphism of $I^*$-algebras.

It is well known that $f^*$ is $*$-morphism. Clearly $f^*I(N)\subset I(M)$. Let $(\phi_n\circ f)_{n\in\NN}$ be a sequence in $f^*I(N)$ converging to an element $\theta$. Now $\theta\in I(M)$, since $I(M)$ is closed, and it defines a function $\psi$ on $f(M)$ such that $\theta=\psi\circ f$. Since $M$ is compact and  both $\theta$ and $f$ are continuous, $\psi$ is continuous, and since $\phi_n\circ f$ converges to $\theta$, $\phi_n$ converges to $\psi$ on $f(M)$, which shows that $\psi$ is isotone. Thus $\psi\in I(f(M))$. Now since $f(M)$ is closed by compacity of $M$, we can use the Tietze-Nachbin theorem to extend $\psi$ to a function $\phi\in I(N)$ such that $\theta=\phi\circ f$. This proves that $f^*I(N)$ is closed.

Now if $J\subset I(N)$ is a pre-isocone of $N$, then $(N,\le_J)$ is a pre-ordered set. If we call $\tilde N$ the quotient of $N$ obtained by identifying the points which are not separated by $J$, we obtain a compact toposet ordered by $p(x)\ll p(y)\Leftrightarrow x\le_J y$, where $p$ is the quotient map. The construction is almost the same as that of proposition \ref{mred} and the verifications are left to the reader. Moreover, since $\le_J$ is finer than the original order on $N$,  $f$ is again a continuous isotony from $M$ to $(N,\le_J)$ and $p$ is a continuous isotony from $(N,\le_J)$ to $(\tilde N,\ll)$. By the same reasoning as above, $(p\circ f)^*I(\tilde N)=f^*p^*I(\tilde N)$ is closed.

Now we just have to show that $p^*I(\tilde N)=J$. 
$$\xymatrix{J&\subset &I(N,\le_J)&\subset&p^*{\mathcal C}(\tilde N)&\subset&{\mathcal C}(N)\cr
\tilde J\ar[u]^{p^*}_{\simeq}&=&I(\tilde N)&\subset&{\mathcal C}(\tilde N)\ar[u]^{p^*}_{\simeq}& }$$

For any function $\phi : N\rightarrow \RR$ going down to $\tilde N$, write $\tilde \phi$ for the unique function such that $\phi=\tilde\phi\circ p$. Then the map $\phi\mapsto \tilde\phi$ is the reciprocal of $p^* : {\mathcal C}(\tilde N)\rightarrow p^*{\mathcal C}(\tilde N)$, which is an isometric $*$-isomorphism. Let $\tilde J$ be the image of $J$ under this map. Then $\tilde J\subset I(\tilde N)$ and satisfies every hypothesis of  corollary \ref{SW2}. Thus we have $\tilde J=I(\tilde N)$. Since $p^*\tilde J=J$ by definition, the result follows.
\end{pf}

We would like to define a functor which goes the other way round. For this we need first to associate a toposet to an $I^*$-algebra.

\begin{defn}
Let $(I,C)$ be an $I^*$-algebra and ${\mathcal S}(C)$ the space of states on $C$. Let $s_1,s_2\in{\mathcal S}(C)$. We define :
\be
s_1\le_I s_2\Leftrightarrow  \forall c\in I,\ s_1(c)\le s_2(c)\label{order}
\ee
\end{defn}

We know from proposition \ref{preord} that $\le_I$  is a preorder. Now suppose $s_1\le_I s_2$ and $s_2\le_I s_1$. Then by definition $s_1$ and $s_2$ coincide on $I$, thus on $\Span(I)$. But $s_1$ and $s_2$ are continuous and then must coincide by axiom 3. Thus $\le_I$ is a partial order on ${\mathcal S}(C)$. In particular $\le_I$ defines a partial order on ${\mathcal PS}(C)$, the space of pure states and on $\chi(C)$, the character space of $C$. Recall that a character of $C$ is $*$-morphism from $C$ to $\CC$ and that $\chi(C)$, ${\mathcal S}(C)$ and ${\mathcal PS}(C)$ are compact for the $*$-weak topology. In the general case one has $\chi(C)\subset {\mathcal PS}(C)$, while there is equality in the commutative case. 

Another space which is associated to $C$ is its structure space $\hat C$, namely the space of equivalence classes of non-trivial irreducible representations modulo unitary equivalence. Through GNS theory, this space can be identified with a quotient of ${\mathcal PS}(C)$. Note that, unless this equivalence relation is identity, which happens exactly when $C$ is commutative, there is no reason $\le_I$ should descend to $\hat C$.

In general we can turn elements of $C$ into continuous functions on ${\chi}(C)$ by the Gelfand transform :
\bea
G(c) &:& C\rightarrow {\mathcal C}({\chi}(C))\cr
& & c\mapsto (m\longmapsto m(c))
\eea

Now it is obvious by definition that $G(c)$ will be isotone when $c\in I$. Thus $G(I)\subset I({\chi}(C))$. But $G(I)$ determines the order $\le_I$ by definition, and this shows that $({\chi}(C),\le_I)$ is a toposet. One can show in the same way that ${\mathcal S}(C)$ and ${\mathcal PS}(C)$ are toposets.

Note that $G$ is always a $*$-morphism. Since $ab$ and $ba$ have the same image, it cannot be into unless $C$ is abelian. If it is the case, then $G$ is a $*$-isomorphism by Gelfand-Naimark theorem. We now state an extension of this theorem to the $I^*$-algebra setting. {\it In the sequel we identify $\hat C$ with $\chi(C)$ when $C$ is commutative.}
 
\begin{thm}\label{GN}
Let $(I,C)$ be an $I^*$-algebra such that $C$ is abelian. Then $\hat C$ is a compact toposet (or equivalently a compact pospace) when it is equipped with the partial order defined by $I$, and the Gelfand transform $G : (I,C)\rightarrow (I(\hat C),{\mathcal C}(\hat C))$ is an isometric isomorphism of $I^*$-algebras.
\end{thm}
\begin{pf}
 By the Gelfand-Naimark theorem we know that $G : C\rightarrow {\mathcal C}(\hat C)$ is an isometric $*$-isomorphism, when $\hat C$ is given the relative weak-$*$ topology, and we have already shown that $\hat C$ is a compact toposet for this topology and the partial order $\le_I$. Now $G$ respects $\wedge$ and $\vee$ since it is a $*$-morphism. Thus $G(I)$ is a subset of $I(\hat C)$ having the following properties : 
\begin{enumerate}
\item It contains the real constants,
\item it is stable by sum, $\wedge$, $\vee$ and mutliplication by a non-negative real,
\item it determines the order. 
\end{enumerate}
In addition to this, $G$ satisfies the third requirement of $I^*$-morphisms since it is a linear isometry. In particular, $G(I)$ is closed. Thus by theorem \ref{SW}, $G(I)=I(\hat C)$. 
\end{pf}

When $A$ is not commutative the Gelfand transform $G : A\rightarrow {\mathcal C}(M)$ is not surjective, where $M={\mathcal PS}(A)$. Would it still be possible that for a given toposet $M$, there exists a noncommutative $I^*$-algebra $(I,A)$ such that $G(I)$ is exactly the set of continuous isotonies on $M$ ? The next proposition answers negatively.

\begin{prop}
Let $A$ be a $C^*$-algebra and $M$ its space of pure states. If there exists an isocone $I$ of $A$ such that $G(I)\subset I(M)$ is stable by $\wedge$ and $\vee$, then $A$ is abelian.
\end{prop}

\begin{pf}
It is clear that $G(I)$ satisfies the hypotheses of the Stone-Nachbin theorem and is closed. Thus $G(I)=I(M)$. Now since $I$ and $I(M)$ satify the third axiom of isocones, we have $G(A)=C(M)$. But this is possible only when $A$ is abelian.  Indeed, $G$ since preserves positivity, it is order-preserving for the order induced by $A_+$ and ${\mathcal C}(M)_+$, respectively. Thus it preserves the lattice operations. From this we can conclude that $A_+$ is a lattice, which proves that $A$ is abelian \cite{sher}.
\end{pf}

Let us now explain in details why the categories ab-${\mathcal IC}$ of abelian $I^*$-algebras and ${\mathcal KT}$ are (dually) equivalent. 

\begin{prop}
Let us call $\widehat{\ }$ the \og mapping\fg\ which associates $(\hat C,\le_I)$ to an abelian $I^*$-algebra $(I,C)$ and the pullback $\hat\phi$ to a morphism $\phi: (I,C)\rightarrow (I',C')$  between two abelian $I^*$-algebras.
Then $\widehat{\ }$ is a contravariant functor from ${\rm ab-}{\mathcal IC}$ to ${\mathcal KT}$. 
\end{prop}
\begin{pf}
The only thing that is not already known from the usual case is that $\hat\phi$ is an isotony. Let us prove this. Take $m,n\in\hat C'$. We have :

\bea
&m\le n&\cr
\Rightarrow &\forall c'\in I',\ m(c')\le n(c')&\mbox{by (\ref{order})}\cr
\Rightarrow &\forall c\in I,\ m(\phi(c))\le n(\phi(c))&\mbox{since }\phi(I)\subset I'\cr
\Rightarrow &\forall c\in I,\ \hat\phi(m)(c)\le \hat\phi(n)(c)\cr
\Rightarrow &\hat\phi(m)\le\hat\phi(n)\nonumber
\eea
\ 
\end{pf}
\underline{Remark} : Note in passing that we nowhere used the fact that $\phi$ satisfies the third requirement of $I^*$-morphisms in the above proof. We will need this fact later.

\begin{prop}
The functors $\widehat{\ }$ and ${\mathcal I}$ realize a dual equivalence between  ab-${\mathcal IC}$ and ${\mathcal KT}$.
\end{prop}
\begin{pf}
We have to show that the named functors are inverse to each other up to natural isomorphisms.

Let us consider a compact toposet $M$ and its image under $\widehat{\ }\circ{\mathcal I}$, i.e. $\widehat{{\mathcal C}(M)}$ with the partial order defined by (\ref{order}). These sets are homeomorphic by the following map :
$$\epsilon_M : M\longrightarrow \widehat{{\mathcal C}(M)}$$
$$x\longmapsto ev_x$$
where $ev_x$ is the evaluation at $x$ : $f\mapsto f(x)$. Now if $x,y\in M$ we have :
\bea
&\epsilon_M(x)\le\epsilon_M(y)&\cr
\Leftrightarrow  &ev_x\le ev_y&\cr
\Leftrightarrow  &\forall f\in I(M),\ f(x)\le f(y)&\mbox{by (\ref{order})}\cr
\Leftrightarrow  &x\le y&\nonumber
\eea
where in the last line we used the fact that $M$ is a toposet.  Thus $M\simeq \widehat{{\mathcal C}(M)}$ as toposets. Moreover, if $f : M\rightarrow M'$ is a continuous isotony between two compact toposets, then $\widehat{f^*} : \widehat{{\mathcal C}(M)}\rightarrow \widehat{{\mathcal C}(M')}$ reads $ev_x\mapsto ev_{f(x)}$ which is clearly a continuous isotony. Thus $\epsilon$ is a natural transformation between the functor $\widehat{\ }\circ {\mathcal I}$ and the identity functor.

Conversely, if $(I,C)$ is an abelian $I^*$-algebra, we have shown in theorem \ref{GN} that $G$ is an isomorphism between $(I,C)$ and $(I(\hat{\mathcal C}),{\mathcal C}(\hat C))$. 
Now we consider :
$${\mathcal I}(\hat\phi) : (I(\hat C),{\mathcal C}(\hat C))\longrightarrow (I(\hat C'),{\mathcal C}(\hat C'))$$
$$f\longmapsto f\circ\hat\phi$$
We already know from the usual case that this commutes with $G$. 
\end{pf}

Thanks to this equivalence, exactly as one thinks of a noncommutative algebra as an algebra of continuous functions on some virtual noncommutative (compact or locally compact) space, one can think of a noncommutative $I^*$-algebra as an algebra of continuous and isotone functions on some virtual noncommutative compact toposet.

We have seen how ideals, quotients and morphisms are defined in $I^*$-algebras. In the commutative case, things are a little simpler. Take two commutative $I^*$-algebras $(I,C)$ and $(I',C')$. Thanks to theorem \ref{GN} we can as well consider the case where $I=I(M)$ and $C=C(M)$ with $M$ a compact toposet, and similarly for the primed objects. A $*$-morphism $\phi$ from $C$ to $C'$ is always the pullback of a continuous $f :M'\rightarrow M$. Now if $\phi(I)\subset I'$ it means that $f$ is a continuous isotony from $M'$ to $M$. In fact, we have already proved this above in showing that $\widehat{\ }$ is a contravariant functor, since, as we have remarked, we nowhere used the third axiom of $I^*$-morphisms. Thus, if $\phi(I)\subset I'$ then $\phi=f^*$ with $f$ a continuous isotony and this proves that $\phi$ is an $I^*$-morphism. Therefore, an $I^*$-morphism between the commutative $I^*$-algebras $(I,C)$ and $(I',C')$ is just a $*$-morphism $\phi$ between $C$ and $C'$ such that $\phi(I)\subset \phi(I')$, and $\phi$ automatically fulfills the third axiom. Consequently, an $I^*$-ideal ${\mathcal N}$ of $(I,C)$ is just any closed ideal of $C$.

Remember that a closed ideal of $C$ is always of the form ${\mathcal N}=\{f\in C|f_{|F}=0\}$, where $F$ is closed subset of $M$, and that the quotient algebra $C/{\mathcal N}$ is canonically isomorphic to ${\mathcal C}(F)$. As expected, this isomorphism carries $I$ to $I(F)$, by the extension property \ref{extension}.

We will now introduce a class of $I^*$-algebra which, as we will show, corresponds to totally ordered toposets.

\begin{defn}
Let $(I,A)$ be an $I^*$-algebra. We will say that $(I,A)$ is {\it minimal} if there exists no $J\subset I$ such that $J$ is an isocone of $A$ and $J\not=I$. In this case, we will also say that $I$ is a minimal isocone of $A$.
\end{defn}

\begin{prop}
Let $M$ be a compact toposet. Then $M$ is totally ordered if, and only if, $(I(M),C(M))$ is a minimal $I^*$-algebra.
\end{prop}
\begin{pf}
Let $M$ be totally ordered and suppose that $(I(M),C(M))$ is not minimal. Then there exists a strict subset $J$ of $I(M)$ such that $J$ is an isocone of $C(M)$. Let $\le$ be the original order on $M$. Then it is clear that for all $x,y\in M$, $x\le y\Rightarrow x\le_J y$. Now if the converse of this implication were true, $J$ would satisfy all conditions of the Stone-Nachbin theorem and we would have $J=I(M)$. Therefore, since $\le$ is a total order, there exist $x,y\in M$ such that $x\le_Jy$ and $y<x$. But $y<x\Rightarrow y<_Jx$, a contradiction.

Conversely, suppose $(I(M),C(M))$ is minimal, and $\le$ is not total. Then there exist $x,y\in M$ such that $x\sim y$ and $x\not=y$. Let $J=\{f\in I(M)|f(x)\le f(y)\}$. Notice first that by definition of a toposet, $x\not=y$ and $x\sim y\Rightarrow \exists g,g'\in I(M)$ such that $g(x)<g(y)$ and $g'(x)>g'(y)$. Thus, there is a $g\in J$ such that $g(x)<g(y)$. We will use this $g$ later in the proof. Since we have a $g'\in I(M)$ such that $g'\notin J$, we also have $J\not=I(M)$.

Now $J$ is obviously a closed and convex cone, and a sublattice of $I(M)$. We just need to show that $\Span(J)$ is dense in $C(M)$. For this we prove that $V=\Span_\RR(J)$ is dense in ${\mathcal C}(M,\RR)$. It is clear that $V=J_+-J_+$ is a subalgebra of ${\mathcal C}(M,\RR)$ containing $\RR 1$, so we only need to show that it separates the points. For this, we take $a\not= b$ in $M$ and $f\in I(M)$ such that $f(a)\not=f(b)$. Now if $f(x)\le f(y)$, or if $g(a)\not=g(b)$ we are done. So we suppose that $f(x)>f(y)$ and $g(a)=g(b)$. Let us define $h\in C(M)$ by 
$$
h=f-2{[ev_x-ev_y](f)\over [ev_x-ev_y](g)}g
$$
We see that $-2{\displaystyle{f(x)-f(y)\over g(x)-g(y)}}>0$, and so $h\in I(M)$. Of course, $h$ being the symmetrical of $f$ with respect to $\ker(ev_x-ev_y)$ in the direction of $g$, we have $h(x)-h(y)<0$ by construction, so $h\in J$. Now we have $h(a)-h(b)=f(a)-f(b)\not=0$. Thus $J$, and then $V$, separate the points.
\end{pf}

\section{Applications}

The first applications of theorem \ref{GN} are to the general (noncommutative) theory of $I^*$-algebras, exactly as in the case of $C^*$-algebras with the Gelfand-Naimark theorem.

\begin{prop}\label{transversecommut}
Let $(I,C)$ be an $I^*$-algebra and $x_1,\ldots,x_p\in I$ such that $[x_i,x_j]=0$ for $1\le i\le j\le p$. Then $(I\cap C^*(x_1,\ldots,x_p),C^*(x_1,\ldots,x_p))$ is an $I^*$-algebra which is isomorphic to $(I(K,\preceq),C(K))$ where $K$ is compact subset of $\RR^p$ and $\preceq$ is a partial order on $K$ which is less fine than the natural product ordering inherited from $\RR^p$.
\end{prop}
\begin{pf}
We already know from $C^*$-algebra theory that there exists such a compact $K\subset \RR^p$ and a $*$-isomorphism $\psi : C^*(x_1,\ldots,x_p)\rightarrow C(K)$ such that $\psi(x_i)$ is the projection on the $i$-th coordinate.

We know that $I\cap C^*(x_1,\ldots,x_p)$ is a pre-isocone. It remains to show that the subalgebra $C^*(x_1,\ldots,x_p)$ is transverse. For this, consider $J=\psi(I\cap C^*(x_1,\ldots,x_p))$. Since $\psi$ is a $*$-isomorphism, $J$ is a pre-isocone of $C(K)$. Let $\preceq$ be the preorder determined by $J$ on $K$. We know from proposition \ref{preord} that $\preceq$ is a partial order if and only if $J$ separates the points. Now this last point is clear since $J$ contains the projections on all coordinates. Thus, $(K,\preceq)$ is a compact toposet, and by corollary \ref{SW2} we have $J=I(K,\preceq)$. This shows that $(I\cap C^*(x_1,\ldots,x_p),C^*(x_1,\ldots,x_p))$ is an $I^*$-algebra which is isomorphic to $(I(K,\preceq),C(K))$. The fact that $\preceq$ is less fine than the product ordering inherited from $\RR^p$ comes from the fact that the projections on coordinates are isotone for $\preceq$.
\end{pf}

If we specialize this proposition to the case $p=1$, we obtain the following.

\begin{cor}\label{propospecident}
Let $(I,C)$ be an $I^*$-algebra and $a\in I$. Then  $(I\cap C^*(a),C^*(a))$ is an $I^*$-algebra. Write $\preceq$ for the order on $\widehat{C^*(a)}$ defined by $I\cap C^*(a)$. Then the map :
$$\theta : (\widehat{C^*(a)},\preceq)\longrightarrow (\sigma(a),\le)$$
$$\phi\longmapsto \phi(a)$$
where $\le$ is the natural order of $\RR$, is an isotone homeomorphism.
\end{cor}

The easy proof of the last part of the corollary is left to the reader. Note that $\theta$ need not be an order-embedding. We will see an example below.

\begin{cor}\label{stabc}
With the same hypotheses, if $f : \sigma(a)\rightarrow \RR$ is continuous and non-decreasing (for the natural order) then $f(a)\in I\cap C^*(a)$.
\end{cor}
\begin{pf}
By theorem \ref{GN} and corollary \ref{propospecident} we have an isomorphism :
$$G : (I\cap C^*(a),C^*(a))\simeq (I(\widehat{C^*(a)},\preceq),\widehat{C^*(a)})$$
Now the map $\theta$ defined above, being an isotone homeomorphism, gives rise to an $I^*$-algebra morphism :
$$\theta^* : (I(\sigma(a),\le),{\mathcal C}(\sigma(a)))\longrightarrow (I(\widehat{C^*(a)},\preceq),\widehat{C^*(a)})$$ which is an isomorphism at the level of $C^*$-algebras. Thus $I(\sigma(a),\le)$ is naturally seen as a subset of $I(\widehat{C^*(a)},\preceq)$. Now remember that we call $f(a)$ the preimage of $\theta^*f$ by the Gelfand-Naimark isomorphism $G$. By the above, $f(a)\in I\cap C^*(a)$.
\end{pf}

Since this may be a little confusing, let us rephrase the proof. By theorem \ref{GN}, $I\cap C^*(a)$ is seen as the set of continuous isotone functions on its spectrum, equipped with the order $\preceq$. By corollary \ref{propospecident}, the spectrum of $C^*(a)$ is identified with the spectrum of $a$ as compact spaces. If we use this identification to transport $\preceq$ we get a partial order that is less fine than the natural order. Thus the continuous functions that are non-decreasing for the natural order will be isotone for $\preceq$ and have a preimage in $I\cap C^*(a)$.

For instance let us take $M=[0;1]$ with the partial order $x\preceq  y$ iff $y=1/2$ or $x=y$. It is clear that $\preceq$ determines a closed subset of $M\times M$ and this shows that $M$ is a compact pospace, hence a toposet. We consider the $I^*$-algebra $(I(M),{\mathcal C}(M))$. Then $I(M)$ is the set of continuous functions on $[0;1]$ reaching their maximum at $x=1/2$. Let $a$ be the function $a(x)=x$ for $x\le 1/2$, and $a(x)=1/2$ for $x\ge 1/2$. It is obvious that $C^*(a)$ is the set of continuous functions on $[0;1]$ which are constant on $[1/2;1]$ and that $I(M)\cap C^*(a)$ contains those which are maximal on $[1/2;1]$. Then $\sigma(a)=[0;1/2]$ and if $f$ is continuous and non-decreasing (for the natural order) on $[0;1/2]$ then $f\circ a=f(a)$ belongs to $I(M)$, but there are other functions in $I(M)\cap C^*(a)$. This shows that $\theta^*$ need not be an isomorphism of $I^*$-algebras, and that $\theta$ need not be an order-isomorphism (which amounts to the same).

\begin{cor}\label{proj}
Let $(I,C)$ be a finite dimensional $I^*$-algebra. Then $I_+$ is the positive span of the set of all the projections  it contains.
\end{cor}
\begin{pf}
Let $a\in I_+$ and $(\lambda_i)_{1\le i\le n}$ its sequence of eigenvalues in increasing order. Then for all $j$ the function $h_j=\sum_{k\ge j}1_{\lambda_k}$ is positive non-decreasing on $\sigma(a)$, so the projection $h_j(a)\in I_+$. We have
$$a=\sum_{k\ge 2}(\lambda_k-\lambda_{k-1})h_k(a)+\lambda_1 h_1(a)$$

Thus $a$ belongs to the positive span of all projections contained in $I_+$.
\end{pf}

Let us see some other applications to proposition \ref{transversecommut}.

\begin{cor}\label{produitcommut}
Let $(I,C)$ be an $I^*$-algebra, and $c,c'\in I_+$. If $[c,c']=0$ then $cc'\in I_+$.
\end{cor}
\begin{pf}
This is obvious by the proposition since for any toposet the set of positive continuous isotonies is stable by product.
\end{pf}

This easily entails the next corollary.

\begin{cor}\label{decadix}
Let $I$ be a pre-isocone of a commutative $C^*$-algebra $A$. Let us write $C^*(I)$ for the smallest sub-$C^*$-algebra of $A$ containing $I$. Then $C^*(I)=\overline{\Span(I)}$.
\end{cor}

This says that when $I$ is a pre-isocone of a commutative $C^*$-algebra $B$, then $(I,C^*(I))$ is an $I^*$-algebra. Remember that given a preordered topological space $M$ we were able to define a reduced space $M_{red}$ which is a toposet, together with a continuous isotone projection $p : M\rightarrow M_{red}$. Now in the case of a pre-ordered $M$, the set $I(M)$ is a pre-isocone of $C(M)$ and it is easy to see that the quotient $M\rightarrow M_{red}$ corresponds to the embedding $(I,C^*(I))\hookrightarrow (I,B)$.

\begin{prop}\label{dualmred}
Let $I$ be a pre-isocone of a commutative $C^*$-algebra $B$ and let $A=C^*(I)$. Then $\hat B_{red}$ is canonically isomorphic to $\hat A$ as a toposet.
\end{prop}
\begin{pf}
Let us just sketch the proof. Thanks to the Gelfand-Naimark theorem we can focus on the case where $B={\mathcal C}(M)$ and $A\simeq {\mathcal C}(N)$, with $N$ a quotient space of $M$. The hypotheses furnish a preorder on $M$ and a partial order on $N$ for which $N$ is a toposet and the quotient map is isotone. Moreover, each element of $I(M)=I$ obviously descend to $I(N)=I$. Therefore, the universal property of the reduced toposet (in the version quoted as a remark after theorem \ref{mred}) is satisfied.
\end{pf}

Another application is the Nachbin compactification. Let $M$ be a completely regular ordered space. That is :

\begin{enumerate}
\item If $\forall x,y\in M,\ \forall f : M\rightarrow [0;1]$ isotone and continuous, $f(x)\le f(y)$ then $x\le y$.
\item $\forall x\in M$, $\forall V$ neighbourhood of $x$, there exist $f$ and $g$ two continuous maps from $M$ to $[0;1]$ such that $f$ is isotone and $g$ is antitone, and such that $f(x)=g(x)=1$, and $f(y)\wedge g(y)=0$ for all $y\in M\setminus V$.
\end{enumerate}
First we remark that the first condition required to be a completely regular ordered space is clearly equivalent to the definition of a toposet. A second obvious remark is that the second condition implies (by considering $h=f\wedge g$) that a completely regular ordered space is completely regular as a topological space, and since it is also Hausdorff by the first condition, it is a Tychonoff space.

Nachbin proved an extension of the Stone-\v Cech compactification for these spaces.

\begin{thm}(Nachbin)

Let $M$ be completely regular ordered space. Then there exists a compact toposet $\nu M$, unique up to isomorphism, such that :
\begin{enumerate}
\item $\nu M$ is an order-compactification of the topological space $M$, i.e. there is a natural inclusion map $i : M\rightarrow \nu M$, which is an order-embedding and a homeomorphism onto its image, and such that $i(M)$ is everywhere dense in $\nu M$.
\item The following universal property is satisfied : for each continuous isotony $F$ from $M$ to a compact toposet $T$, there exists a continuous and isotone $\tilde F : \nu M\rightarrow T$ such that $F=\tilde F\circ i$ (i.e. the following diagram commutes).
\end{enumerate}
$$\xymatrix{M\ar@{->}[r]^F \ar@{->}[d]_i&T\cr \nu M\ar@{-->}[ur]_{\tilde F}&}$$
\end{thm}

We wish to recover this theorem by an isocone construction. We could consider the couple $(I_b(M),B)$, where $B={\mathcal C}_b(M)$ stands for the set of continuous bounded functions from $M$ to $\CC$, but it is not necessarily an $I^*$-algebra. Indeed, let us call $A=\overline{\Span(I_b(M))}$. Then $A$ is a closed subalgebra of $B$, and hence a $C^*$-algebra itself. But $A\not=B$ in general. For instance if $M=\RR$ with the natural order, then it is easy to see that $f\in A\Rightarrow \lim_{x\rightarrow \pm\infty}f(x)$ both exist. On the other hand, $I_b(M)$ is clearly a pre-isocone. Thus, $(I_b(M),A)$ is an $I^*$-algebra by corollary \ref{decadix}. By theorem \ref{GN} we have $(I_b(M),A)\approx (I(\hat A),{\mathcal C}(\hat A))$ for a compact toposet $\hat A$. What is the relation between $M$ and $\hat A$ ? Using proposition \ref{dualmred} we clearly have $\hat A=\hat B_{red}$. Moreover we know that $\hat B=\beta M$ is the Stone-\v Cech compactification of $M$. We will show below that $\hat A$ is also the Nachbin compactification of $M$.

\begin{thm}
Let $M$ be a toposet which is completely regular as an ordered space. Let $A$ and $B$ be defined as above. Write $\beta M$ for the Stone-\v Cech compactification $\hat B$ of $M$, and $\nu M$ for its Nachbin compactification. Then :
\begin{enumerate}
\item $\hat A$ is canonically isomorphic to $\nu M$.
\item Let $\preceq$ be the preorder on $\beta M$ defined by $I_b(M)$, i.e. $\phi\preceq\phi'$ if, and only if, $\forall a\in I_b(M),\ \phi(a)\le \phi'(a)$. Then $\nu M$ is the reduced toposet $\beta M_{red}$ with respect to this preorder.
\end{enumerate}
\end{thm}
\begin{pf}
The second assertion is obvious given the first. We summarize what we already know from usual $C^*$-theory and from what we said above in the following commutative diagram, where $\hat A$ is a compact toposet, $p$ a continuous isotone surjection, the continuous maps $\epsilon$ and $i$ are given by $x\mapsto ev_x$ and $\epsilon$ is injective.

$$\xymatrix{ &&  \beta M\ar@{->>}[dd]^p\cr 
&&\cr
M\ar@{^{(}->}[uurr]^\epsilon \ar@{->}[rr]^i  && \hat A=\beta M_{red}}$$

Now $i$ is also injective precisely because $M$ is a toposet. If $f\in{\mathcal C}(\hat A)\simeq A$ vanishes on $i(M)$, then clearly $f=0$. This shows that $i(M)$ is dense in $\hat A$, since we can use Urysohn's lemma in the compact and thus completely regular space $\hat A$.

We now show that $i$ is a homeomorphism onto its image. Indeed, take $U$ an open set in $M$ and $x\in M$. Then use complete regularity of the ordered space $M$ to find $f$ and $g$ as in the definition (note that $g\in A$). The set $O=\{\phi\in\hat A|\phi(f)>0,\phi(g)>0\}$ is a $*$-weak open set. This shows that $V=i(M)\cap O$ is an open set for the induced topology on $i(M)$ and furthermore, it is clear that $i(x)\in V$ and thanks to properties of $f$ and $g$, $V\subset i(U)$. Thus, $i(U)$ is open, proving that $i$ is a homeomorphism.

Now the partial order on $\hat A$ is defined by :

$$\phi\le \phi'\Leftrightarrow \forall f\in I,\ \phi(f)\le \phi'(f)$$
Specializing to $\phi=ev_x$ and $\phi'=ev_y$, we see that $i$ is an order-embedding. This already shows that $\hat A$ is an ordered compactification of $M$.

Now let us prove that $\beta M_{red}$ satisfies the same universal property as $\nu M$. For this, let $F$ be a continuous isotony to a compact toposet $T$. Using the fact that $\Span(I(T))$ is dense in ${\mathcal C}(T)$, we see that for any $\psi\in{\mathcal C}(T)$, $\psi\circ F\in A$. Take $\chi\in\hat A$ and define $\tilde F(\chi):=[\psi\mapsto\chi(\psi\circ F)]$. It is clear that $\tilde F$ is continuous map from $\hat A$ to $\widehat{{\mathcal C}(T)}=T$ extending $F$, where we allow ourselves to identify spaces isomorphic in a well-behaved way with respect to every map in sight. Now let $\chi\le\chi'$ in $\hat A$. Then $\tilde F(\chi)\le\tilde F(\chi')$ iff $\forall \phi\in I(T)$, $\chi(\phi\circ F)\le\chi'(\phi\circ F)$, which is true by definition of the order on $\hat A$. Thus $\tilde F$ is isotone and the theorem is true.
\end{pf}

\section{Noncommutative causal toposets}
We now look for a noncommutative generalization of causal sets. Let us first recall the definition of these objects.

\begin{defn}
A causal set is a poset $P$ which is locally finite, i.e. $\forall x,y\in P$, $\{z\in P| x\le z\le y\}$ is finite. We will say that $P$ is a causal toposet if $P$ is a locally finite toposet.
\end{defn}

Note that a causet does not have a topology defined on it. However, if we stick to finite causets, the discrete topology naturally endows them with a toposet structure. The dual object of a finite causet is obviously a finite dimensional commutative $I^*$-algebra. Thus, the noncommutative generalization of a finite causet is just a finite dimensional $I^*$-algebra. This is a special case, but it is the most useful one. Indeed, a causet appearing in a quantum theory of gravity will certainly be finite, or at most be a countable inductive limit of finite causets. In this last case, the noncommutative version would be a countable projective limit of finite dimensional $I^*$-algebras. However, this definition is unsatisfactory because the inductive limit is not canonically defined. For this reason, and also for the sake of generality, we need to look for a dual formulation of the local finiteness condition of causal sets.

First, we recall a few definitions. In a poset $P$, the set $\{z\in P| x\le z\le y\}$ is called a {\it closed interval}. A poset is called {\it bounded} if it has a greatest and a lowest element. A subposet $M\subset P$ is called a {\it gem}\footnote{This terminology is particular to us, and has been suggested by the shape of the Hasse diagram of such a subposet. We hope not to be in conflict with any existing vocabulary.} if it is a bounded poset for the induced order. Note that a subposet $M$ is called bounded if it has a upper and a lower bound in $P$, so that a gem is not just a bounded subposet. Thus we can say that a causal set is a poset in which every closed interval is finite. Since a closed interval is a gem, and each gem is contained in a closed interval, it is clear that $P$ is a causal set iff $P$ is a poset in which every gem is finite.

\begin{prop}
If $M$ is a compact bounded toposet\footnote{That is, $M$ is compact as a topological space and bounded as a poset.}, then :
\begin{enumerate}
\item $\forall f,g\in I_+(M)$, $\|f+g\|=\|f\|+\|g\|$
\item $\forall f,g\in I_+^\times(M)$, $\|f^{-1}+g^{-1}\|=\|f^{-1}\|+\|g^{-1}\|$
\end{enumerate}
Conversely, if $(I,C)$ is a commutative $I^*$-algebra such that :
\begin{enumerate}
\item $\forall f,g\in I_+$, $\|f+g\|=\|f\|+\|g\|$
\item $\forall f,g\in I_+^\times$, $\|f^{-1}+g^{-1}\|=\|f^{-1}\|+\|g^{-1}\|$
\end{enumerate}

then $(\hat C,\le_I)$ is a compact bounded toposet.
\end{prop}
\begin{pf}
For the first part, let $a$ be the least and $b$ the greatest element of $M$. Then every element of $I(M)$ will attain its minimum at $a$ and its maximum at $b$. The result follows.

For the converse, we can use theorem \ref{GN} and consider only the case where $C=C(M)$, $I=I(M)$ with $M$ a compact toposet. Let us show that $M$ is bounded. Let $f,g\in I_+$. Then $\|f+g\|$ is attained at some $x\in M$ by compacity. So, using the hypothesis, $\|f+g\|=f(x)+g(x)=\sup f+\sup g$ which easily entails that $f$ and $g$ reach their respective suprema at the same point $x$. By induction, one can show that $\forall f_1,\ldots,f_n\in I_+$, $\exists x\in M$ such that $f_1,\ldots,f_n$ all reach their sup at $x$. For all $f\in C_+$, write $S(f)=\{x\in M|f(x)=\|f\|\}$. Obviously, $S(f)$ is a closed subset of $M$. We have shown that 
$$\forall f_1,\ldots,f_n\in I_+, \bigcap_{i=1}^nS(f_i)\not=\emptyset$$
therefore
$$\bigcap_{f\in I_+}\ S(f)\not=\emptyset$$
by compacity of $M$. Let $b$ belong to this intersection. Let $x\in M$. For all $f\in I$, $f+\|f\|\in I_+$ and one has $f(x)+\|f\|\le f(b)+\|f\|$ by definition of $b$. Since $M$ is a toposet, this shows that $x\le b$. Consequently, $b$ is a maximal element for $M$.

The proof that there exists a minimal element is similar, using $\sup(1/f)=1/\inf(f)$, which is valid for $f\in I_+^\times$.
\end{pf}
 
Since the dual formulation of boundedness does not use commutativity of the algebra, this motivates the following definition :

\begin{defn}
An $I^*$-algebra $(I,A)$ is said to be co-bounded if
\begin{enumerate}
\item $\forall a,b\in I_+$, $\|a+b\|=\|a\|+\|b\|$
\item $\forall a,b\in I_+^\times$, $\|a^{-1}+b^{-1}\|=\|a^{-1}\|+\|b^{-1}\|$
\end{enumerate}
\end{defn}

However it is much easier in practice to use the following equivalent formulation in terms of pure states :

\begin{enumerate}
\item $\forall a,b\in I_+$, ${\displaystyle\sup_{\phi\in{\mathcal PS}(A)}\phi(a+b)=\sup_{\phi\in{\mathcal PS}(A)}\phi(a)+\sup_{\phi\in{\mathcal PS}(A)}\phi(b)}$,
\item $\forall a,b\in I_+$, ${\displaystyle\inf_{\phi\in{\mathcal PS}(A)}\phi(a+b)=\inf_{\phi\in{\mathcal PS}(A)}\phi(a)+\inf_{\phi\in{\mathcal PS}(A)}\phi(b)}$.
\end{enumerate}

For a noncommutative example of a co-bounded $I^*$-algebra, take any $I^*$-algebra $(I,C)$ and then consider the $C^*$-algebra $A=\CC\oplus C\oplus \CC$, and define $F$ to be set of hermitian elements of $A$ of the form $\lambda\oplus a\oplus \mu$, with $a\in I$, $\lambda\ge\|a\|$ and $\mu\le{\displaystyle\inf_{\phi\in{\mathcal PS}(A)}\phi(a)}$. It is readily verified that $(F,A)$ is a co-bounded $I^*$-algebra. In the commutative case this construction is dual to the adjunction of topologically isolated top and a bottom elements to a compact toposet. 

Now every inclusion $M\hookrightarrow P$ of a gem into a toposet corresponds to surjective $I^*$-morphism going in the other direction. Thus the dual notion of co-gem is naturally defined as follows.

\begin{defn}
Let $(I,A)$ be an $I^*$-algebra. A co-gem of $(I,A)$ is a quotient of $(I,A)$ which is co-bounded.
\end{defn}

And we finally arrive at the noncommutative version of a \og causal toposet\fg, that is, a toposet which is also a causet.

\begin{defn}
A causal $I^*$-algebra is an $I^*$-algebra such that each of its co-gems is finite-dimensional. 
\end{defn}

To go from causal toposets to causal sets, we could try to see if a causet is canonically a toposet for a topology generated by the order in some way. This is a direction we look forward to exploring. Note that the correspondence between causal toposets and causal algebras is functorial, contrarily to the correspondence with incidence algebras explored in \cite{sork2}. 

\section{Isocones on $M_2(\CC)$}

We will now focus on the simplest possible cases of noncommutative $I^*$-algebras, namely ones with $M_2(\CC)$ as algebra. For a finite-dimensional abelian $C^*$-algebra, there is only a finite number of possible isocones, in one-to-one correspondence with the different partial orders on the character space, which is finite. Now we will see already with the case of $M_2(\CC)$ that there is an (uncountably) infinite number of different isocones.

We use the basis of $\reel M_2(\CC)$ given by the Pauli matrices :

$$\sigma^0=I_2=\pmatrix{1&0\cr 0&1},\ \sigma^1=\pmatrix{0&1\cr 1&0},\ \sigma^2=\pmatrix{0&-i\cr i&0},\ \sigma^3=\pmatrix{1&0\cr 0&-1}$$

We let $T=\Span_\RR\{\sigma^1,\sigma^2,\sigma^3\}$ be the space of traceless matrices. If $I$ is an isocone of $M_2(\CC)$ (it makes no difference if it is weak or strong, as we will see), then $C=I\cap T$ is obviously a closed convex cone in $T$. Moreover, since $I$ satisfies axiom \ref{dense}, it must contain a $4$-ball $B$. If $a$ is the center of $B$, then $(B-{1\over 2}\tr(a)\sigma^0)\cap T$ is a $3$-ball contained in $C$. Thus we have a well-defined map $\phi : {\mathcal I}\rightarrow \Gamma$, from the set ${\mathcal I}$ of (strong or weak) isocones of $M_2(\CC)$ to the set $\Gamma$ of closed convex cones with non-empty interior in $T\simeq \RR^3$.

\begin{thm}\label{classisom2c}
The map $\phi$ is a bijection with $\phi^{-1}$ given by $C\mapsto C+\RR\sigma^0$.
\end{thm}
\begin{pf}
The fact that $C\mapsto C+\RR\sigma^0$ defines an inverse to $\phi$ is immediate provided we prove that $C+\RR\sigma^0$ is an isocone when $C$ belongs to ${\mathcal I}$. All axioms except \ref{meetandjoin} are directly seen to be satisfied. It turns out that \ref{meetandjoin} is then automatically true. Indeed, it is an easy exercise to prove that for all $a,b\in \reel M_2(\CC)$ : 
$$a\vee b=\alpha a+(1-\alpha) b+\beta \sigma^0$$
and
$$a\wedge b=(1-\alpha) a+\alpha b-\beta \sigma^0$$
with $\alpha\in[0;1]$ and $\beta\in\RR_+$. This shows that for a subset of $\reel M_2(\CC)$, the axioms \ref{somme}, \ref{scal} and \ref{constants} of isocones imply  \ref{meetandjoin}. In particular, there is no difference between strong and weak isocones in $M_2(\CC)$.
\end{pf}

Since a convex cone in $T$ is characterized by its intersection with the unit sphere $\Sigma$ of this space, it is also possible to classify the isocones of $M_2(\CC)$ by the subsets of $\Sigma$ which are closed, have a non-empty interior (for the topology of the sphere), and are either the whole sphere or a geodesically convex subset of a closed half-sphere. If $K$ is a subset of $\Sigma$ of the sort just mentioned, we write $I_K$ for the isocone it generates, that is : $I_K=\RR_+K+\RR\sigma^0$.

It is is easy to classify $I^*$-algebra structures up to isomorphism. Indeed, any $*$-automorphism of $M_2(\CC)$ is of the form $a\mapsto pap^*$, with $p\in SU(2)$ and induce a rotation of $T$. Thus the list of isocones given by theorem \ref{classisom2c}, up to rotations, gives the list of isomorphism classes of $I^*$-algebra structures on $M_2(\CC)$. In particular, the automorphisms of $(I_K,M_2(\CC))$ are given by those rotations which leave $K$ invariant. In the commutative case, the automorphisms correspond by duality to the order-automorphisms of the underlying toposet. We see that the automorphism group is much larger in the noncommutative case.

Given a set ${K}$ of the above kind, and a corresponding isocone $I_{K}$, one can wonder what is the induced partial order relation $\le_{I_{K}}$ on the set of pure states of $M_2(\CC)$. It turns out that it has a simple geometric interpretation.

Recall that the pure states space of $M_2(\CC)$ is $P^1(\CC)$, and that for each unit vector $\xi$ in $\CC^2$ with class $[\xi]$ in $P^1(\CC)$ the corresponding state is given by $\phi_{[\xi]}(m)=\xi^*m\xi$. Now it is easy to see that 
$$\phi_{[\xi]} \le_{I_{K}}\phi_{[\eta]}\Leftrightarrow \forall m\in{K},\ \phi_{[\xi]}(m)\le\phi_{[\eta]}(m)$$

This shows that  $\le_{I_{K}}=\le_{K}$. Moreover, if we write $\vect{m}$ for the column vector of coordinates of $m$ in the basis $(\sigma^1,\sigma^2,\sigma^3)$, then by a direct calculation we find that $\phi_{[\xi]}(m)=\vect{m}.\vect{h}([\xi])$, where $\vect{h} : P^1(\CC)\rightarrow S^2$ is the isomorphism between $P^1(\CC)$ and $S^2$ induced by the Hopf fibration :

$$\vect{h}([\xi_1,\xi_2])=\pmatrix{2\reel(\bar\xi_1\xi_2)\cr 2\Im(\bar\xi_1\xi_2)\cr |\xi_1|^2-|\xi_2|^2)}$$

Thus we have $\phi_{[\xi]}\le_{K}\phi_{[\eta]}\Leftrightarrow \forall m\in{K},\ \vect{m}.\vect{h}([\xi])\le\vect{m}.\vect{h}([\eta])$. If we set $d(a,b)$ for the geodesic distance between the points $a$ and $b$ in $S^2$, we then have 

\be
\phi_{[\xi]}\le_{K}\phi_{[\eta]}\Leftrightarrow \forall{m}\in{K},\ d(\vect{m},\vect{h}([\xi]))\ge d(\vect{m},\vect{h}([\eta]))
\ee

Now, since the Hopf fibration takes the canonical metric of the sphere to the Fubini-Study metric of $P^1(\CC)$ (up to an irrelevant scale factor), we finally obtain

\be
\phi_{[\xi]}\le_{K}\phi_{[\eta]}\Leftrightarrow \forall{x}\in{K'},\ \delta(x,[\xi])\ge \delta(x,[\eta])
\ee

where $K'=\vect{h}^{-1}(K)$ and $\delta$ is the Fubini-Study metric.

\underline{Remark} : This kind of order relation can be generalized to any metric space. Let us say that a subset $K$ of a metric space $E$ is GPS-complete if any point $x$ in $E$ is uniquely characterized by the set of distances $d(x,z)$ with $z\in K$. Then, define the GPS-ordering $\le_K$ to be the relation $x\le_K y$ iff $d(x,z)\le d(y,z)$ for all $z\in K$. This relation is an order relation thanks to the GPS-completeness property of $K$, and it turns $E$ into a toposet since the functions $x\mapsto d(x,z)$ are continuous, isotone, and determine the order. This furnishes a abundant supply of toposets.

Another question we can ask is what is the order induced on the spectrum of a transverse normal element $n$. Take $n$ such that $n^*n=nn^*$, and call $p$ a unitary matrix diagonalizing $n$ such that $n=pdp^*$ with $d=\pmatrix{\lambda_1&0\cr 0&\lambda_2}$ with $\lambda_1\not=\lambda_2$. Let us write $\xi=\pmatrix{z_1\cr z_2}$ and $\xi'=\pmatrix{-\bar z_2\cr \bar z_1}$ for the columns of $p$. With this choice one has $p\sigma^3p^*=\vect{h}([\xi])=-\vect{h}([\xi'])$. Thus we have

\be
x={1\over 2}((\lambda_1+\lambda_2)\sigma^0+(\lambda_1-\lambda_2)\vect{h}([\xi]))\label{formula}
\ee

and $\reel C^*(n)=\Span_\RR\{\sigma_0,\vect{h}([\xi])\}$. We see that $n$ is transverse if and only if the line $\RR \vect{h}([\xi])$ cuts ${K}$. Since the intersection point will act as an isotone function on $\sigma(n)$, and $\vect{h}([\xi])=p\sigma^3p^*$ acts on $\sigma(n)$ by $\lambda_1\mapsto 1$ and $\lambda_2\mapsto -1$, we have, writing $\preceq$ for the induced order on $\sigma(n)$ in case $n$ is transverse 
:
\begin{enumerate}
\item If $\vect{h}([\xi])\in {K}$, then $\lambda_2\preceq \lambda_1$.
\item If $-\vect{h}([\xi])\in {K}$, then $\lambda_1\preceq\lambda_2$.
\item If both belong to ${K}$ then $\lambda_1\sim \lambda_2$.
\item If neither belong to ${K}$ then $n$ is not transverse.
\end{enumerate}

In particular, for a hermitian element $x$, we see thanks to (\ref{formula}) that $x$ is transverse if and only if $x\in I$, or $-x\in I$. If $x\in I$ and $-x\notin I$ then we can directly check that the order induced on $\sigma(x)=\{\lambda_1;\lambda_2\}$ is the natural order of $\RR$, in accordance with proposition \ref{propospecident}, whereas it is the dual order if $-x\in I$ and $x\notin I$.

We can give a quantum mechanical interpretation of this fact. Let us consider a two-state system characterized by the values of an observable $x$. When observed, the system can be found in the state $[t]$ (for top), corresponding to the largest eigenvalue of $x$, or $[b]$ (for bottom), corresponding to the smallest eigenvalue. We can see $[t]$ and $[b]$ as the north and south poles of the pure state space $P^1(\CC)$ identified with a sphere. Now any isocone of the form $I_K$ with $K$ containing the north pole (that is, $\vect{h}([t])$) in its interior will induce a partial order on $P^1(\CC)$ which will extend the relation $[b]<[t]$ to the quantum superpositions of these two states. A natural choice for $K$ would be the geodesic disk of radius $\epsilon$, for any $\epsilon\in]0;\pi/4]$ (remember that $\pi/4$ is the distance from the pole to the equator in the Fubini-Study metric). It is then easy to check that, for any pure state $[\xi]$, if $d([b],[\xi])<2\epsilon$ then $[b]\sim [\xi]$, if $d([b],[\xi])\ge 2\epsilon$, then $[b]\le_K [\xi]$, and similarly, if $d([t],[\xi])<2\epsilon$, then $[t]\sim[\xi]$, and if $d([t],[\xi])\ge 2\epsilon$, $[\xi]\le_K[t]$. Now when the system is in the state $[\xi]$, the probability that it will collapse to $[t]$ when $x$ is measured is $\cos^2(d([\xi],[t]))$. Thus, we have just seen that for the order $\le_K$ with our choice of $K$, a state $[\xi]$ is smaller that $[t]$ (or larger than $[b]$) if the probability that it will collapse to the other eigenstate is larger than some level dependending on $\epsilon$, namely $1-\cos^2(2\epsilon)$.

\section{Conclusion}

In this paper we have seen that it is possible to add a line corresponding to a class of ordered space in the algebra/geometry dictionary set up by the noncommutative geometry programme. Not only is it possible, but it is indeed quite natural, at least at the level of generality where we have settled. 

On the mathematical side many questions arise. To quote only a few, we would like to know if the nice geometric characterization of isocones on $M_2(\CC)$ survives in any way in more elaborate examples. It is not completely hopeless. Indeed, though the pure state space ${\mathcal PS}(A)$ of a $C^*$-algebra $A$ does not suffice to characterize the algebra in the noncommutative case, it is possible to add more structures on ${\mathcal PS}(A)$ (namely its uniform structure, an orientation and transition probabilities, see \cite{schultz}) so that it does. In our context, it is natural to expect that $(I,A)$ will be completely determined by its pure states space if we add a toposet structure compatible in some way with these data.  It would also be interesting to know if the definition of morphisms can be relaxed. It is also clear that more examples of noncommutative $I^*$-algebras are needed. Work in this direction is in progress.

\section{Acknowledgement}

We would like to thank the anonymous referee for his interesting suggestions and remarks.

\end{document}